\documentclass[
 aip,
rsi,
 amsmath,amssymb,
 reprint,
nofootinbib,
floatfix,
]{revtex4-1}

\usepackage{url}

\usepackage{graphicx}
\usepackage{dcolumn}
\usepackage{bm}
\usepackage{color}
\usepackage[dvipsnames]{xcolor}
\usepackage{soul}

\usepackage[caption=false]{subfig}

\usepackage{hhline}

\usepackage{verbatim}

\begin{document}

\preprint{UWA/001-DMW}

\title[Consistency in Echo-State Networks]{Consistency in Echo-State Networks}

\author{Thomas Lymburn}
\email{thomas.lymburn@research.uwa.edu.au}
\affiliation{Complex Systems Group, Department of Mathematics and Statistics, Faculty of Engineering and Mathematical Sciences, The University of Western Australia, Crawley, Western Australia 6009, Australia}

\author{Alexander Khor}
\affiliation{Complex Systems Group, Department of Mathematics and Statistics, Faculty of Engineering and Mathematical Sciences, The University of Western Australia, Crawley, Western Australia 6009, Australia}

\author{Thomas Stemler}
\affiliation{Complex Systems Group, Department of Mathematics and Statistics, Faculty of Engineering and Mathematical Sciences, The University of Western Australia, Crawley, Western Australia 6009, Australia}

\author{D{\'e}bora C. Corr{\^e}a}
\affiliation{Complex Systems Group, Department of Mathematics and Statistics, Faculty of Engineering and Mathematical Sciences, The University of Western Australia, Crawley, Western Australia 6009, Australia}

\author{Michael Small}
\affiliation{Complex Systems Group, Department of Mathematics and Statistics, Faculty of Engineering and Mathematical Sciences, The University of Western Australia, Crawley, Western Australia 6009, Australia}
\affiliation{Mineral Resources, CSIRO, Kensington, Western Australia 6151, Australia}

\author{Thomas J{\"u}ngling}
\affiliation{Complex Systems Group, Department of Mathematics and Statistics, Faculty of Engineering and Mathematical Sciences, The University of Western Australia, Crawley, Western Australia 6009, Australia}

\date{\today}

\begin{abstract}
Consistency is an extension to generalized synchronization which quantifies the degree of functional dependency of a driven nonlinear system to its input.
We apply this concept to echo-state networks, which are an artificial-neural network version of reservoir computing.
Through a replica test we measure the consistency levels of the high-dimensional response, yielding a comprehensive portrait of the echo-state property.
\end{abstract}

\maketitle


\begin{quotation}
When a nonlinear dynamical system is externally modulated by an information-carrying signal, its erratic response hides an intricate property: Consistency.
It is difficult to estimate from time series whether or not the variability in the output is entirely determined by the driving signal.
For autonomous chaotic systems it is well-known that their inherent instability gives rise to a certain level of unpredictability.
For a driven system, this means that a part of the variability of its output does not depend on the drive.
Consistency quantifies the degree of this dependency through a replica test.
The nonlinear system is repeatedly driven by the same signal, and the corresponding responses are compared.
We apply this concept to echo-state networks, a class of artificial neural networks with a fixed random internal connectivity.
Such networks have been successfully utilized for sequential processing tasks like nonlinear time series prediction and spoken digit recognition.
Studying the consistency property allows for a more comprehensive understanding of the dynamical response and for tailoring the network systematically towards enhanced functionality and a wider range of applications.
\end{quotation}


\section{Introduction}

Synchronization is a common phenomenon in interacting nonlinear oscillators that has been studied for almost three decades~\cite{Boccaletti:02,Arenas:08}.
The mutual or directed interaction can lead to several forms of entrainment of the trajectories.
Different degrees of relationships have been discussed, like complete synchronization (CS) and phase synchronization (PS)~\cite{Pikovsky:01}.
For generalized synchronization (GS)~\cite{Rulkov:95,Moskalenko:12}, however, only the presence and form of a functional relationship have been analyzed, but no corresponding weaker form of synchronization received any considerable attention~\cite{Giacomelli2010a}.
This least-studied case may be the most prevalent in natural systems and also of high relevance for novel forms of neuro-inspired computation.

When a nonlinear system is driven by an external signal - like neurons which are excited by real-world stimuli - the nature of the dependency between drive and response is a very important and challenging aspect.
In the field of \emph{reservoir computing} (RC), which overlaps with recurrent neural networks (RNN), dynamical systems are employed to tailor functions on sequential data~\cite{Lukosevicius:09,Appeltant:11,Konkoli:18}.
In contrast to feedforward structures, which constitute the majority of present artificial neural networks (ANN), dynamical systems are known to develop instabilities.
This property is generally undesired and hard to control, which is one of the reasons for the marginal existence of RC and RNN.
In this work, we connect a concept from nonlinear science, namely consistency, and the associated replica scheme to echo-state networks (ESN)~\cite{Jaeger:01,Jaeger:04}, which are a particular flavor of RC that is based on RNN.

The concept of consistency emerged mainly within the last decade as an approach to introduce nonlinear science methodology to a broader domain in which the response of a nonlinear dynamical system to arbitrary signals plays a crucial role~\cite{Uchida:04,Uchida2008a,Kanno2012,Oliver:15,Nakayama:16,Bueno:17,Jungling:18}.
Consistency is based on the replica test, in which a nonlinear dynamical system is repeatedly driven with the same signal.
The test is an adaption of the Abarbanel test for GS~\cite{Abarbanel:96}, in which at least two identical units are driven simultaneously.
In each version of the replica test, the different responses, which in theory just differ in their initial conditions, are compared, typically by means of a correlation coefficient which measures the degree of consistency~\cite{Uchida:04,Jungling:18}.
The term consistency has not been defined rigorously yet, and there are currently different possible interpretations.
We elaborate on this issue in Sec.~\ref{sec:consistency}.
A major goal of this work is to contribute to an enhanced understanding of consistency, in particular for large dynamical systems such as the ESN.

Echo-state networks are a computationally feasible RC paradigm, which is distinguished from other RNNs by a simplified training procedure~\cite{Jaeger:01,Lukosevicius:09}.
The main idea of RC is to utilize the response of a large dynamical system, the reservoir, to generate nonlinear features in a high-dimensional space.
Reservoir computing generally employs \emph{physical} dynamical systems and thus belongs to the field of unconventional computation (UC)~\cite{Konkoli:18}.
The reservoir in ESN is a numerical model, typically given by a realization of a random network of dynamical nodes with sigmoid activation functions.
Without driving signal the reservoir is typically designed to reside in a stable steady state, and the transient activation during injection of the input is recorded.
This random nonlinear embedding of the signal in a large dimension facilitates certain regression or classification tasks \cite{Pathak2017,Lu2017,Zimmermann2018}.

Consistency is a property of the nonlinear response of ESN which is related to conditional stability.
Stability in a driven nonlinear system typically refers to the spectrum of conditional Lyapunov exponents (CLE) and the corresponding Lyapunov vectors, most importantly the largest exponent~\cite{Pecora:90,Pyragas:97}.
All these quantities depend on both properties of the driven dynamical system and the driving signal~\cite{Heiligenthal:11}.
The CLE may often appear to be tightly linked to consistency, however, both are complementary characteristics of the driven system~\cite{Oliver:15}.
In ESN, conditional stability is a synonym for the \emph{echo-state property}, which is a well-known central characteristic that is often referred to as a necessary criterion for the function of the network as a reservoir.
Under ideal conditions, the echo-state property is equivalent to \emph{complete consistency}, which refers to identical responses to repetitions of the driving signals~\cite{Oliver:15,Jungling:18}.
A few attempts to obtain a deeper insight into the echo-state property have been presented~\cite{Verstraeten:09,Yildiz:12,Manjunath:13}, as well as a mean-field theory for the signal propagation~\cite{Massar:13}, and the use of reservoirs in non-stationary regimes~\cite{Marquez:18}. 

In this work, we employ the ESN in a parameter regime beyond the typical ranges that guarantee the echo-state property, as well as with intrinsic noise which has a similar effect.
The replica test is implemented for realizations of Gaussian white noise as a scalar driving signal.
We distinguish between the micro-level consistency of the reservoir nodes and the emergent consistency at the readout level, thus yielding a comprehensive portrait of the consistency property.
This perspective is important for the response characterization in neuronal microcircuits, in which many factors lead to a noisy micro-level but still allow for a robust functionality, see also the concept of coarse coding~\cite{MacLennan:18,Sanger:96,Rumelhart:86}.
Moreover, our approach is applicable to neuro-inspired technical applications in which an experimental access to the reliability of the response systems is required.
We elaborate on the general consistency property in Sec.~\ref{sec:consistency}.
In Section~\ref{sec:memory}, we investigate the relationship between consistency and the fading memory in an ESN.
We finally introduce a \emph{consistency profile} in Sec.~\ref{sec:profile} which demonstrates how an injected signal propagates and fades in the fluctuating neural medium.


\section{The consistency property}
\label{sec:consistency}

The replica test and the consistency measure have so far been applied only to low-dimensional systems and scalar time series.
When transferring the concept to complex dynamical systems like an ESN, one encounters several new aspects.
We approach these with a general driven system acting as a reservoir
\begin{equation}
\dot{\mathbf{x}}(t)=\mathbf{f}(\mathbf{x}(t),\mathbf{u}(t),\mathbf{q})\;.
\label{eq:reservoir}
\end{equation}
Here, $\mathbf{x}(t)\in\mathbb{R}^N$ is the state of the reservoir $\mathcal{X}$ through continuous time, which we will also refer to as a network of nodes representing the degrees of freedom~\cite{Appeltant:11,Grigoryeva:15}.
The vector $\mathbf{u}(t)\in\mathbb{R}^L$ is the multivariate driving signal, and $\mathbf{q}\in\mathbb{R}^M$ is a set of fixed parameters which control internal wiring of the reservoir, the input injection, as well as the shape of the nonlinearity.

\begin{figure}
\begin{center}
\includegraphics[width=\columnwidth]{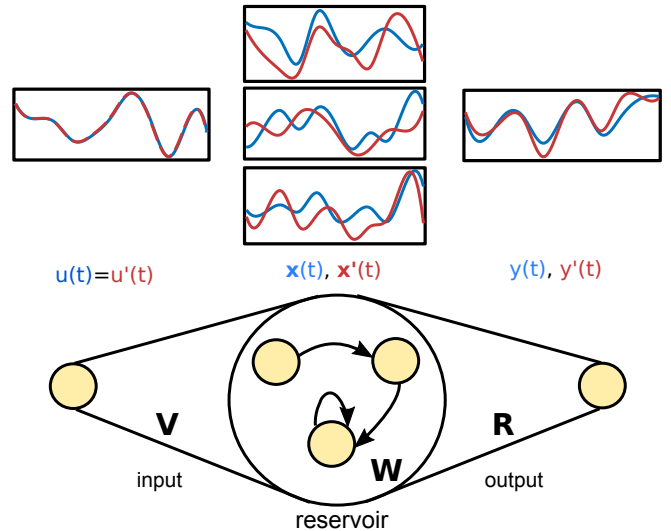}
\end{center}
\caption{Illustration of the replica scheme for a reservoir computer according to Eqs.~\eqref{eq:reservoir}--\eqref{eq:GammaR}. Bottom: Topology with input node, reservoir, and output node and connectivities $\mathbf{V},W,\mathbf{R}$. Top: Original (blue) and repetition (red) of signals at corresponding nodes.}
\label{fig:replica}
\end{figure}

In a basic replica test, an identical copy $\mathcal{X}'$ of the reservoir $\mathcal{X}$ is simultaneously driven with the same signal $\mathbf{u}(t)$, but starting from different initial conditions $\mathbf{x}'(t_0)\neq\mathbf{x}(t_0)$.
Alternatively, the same reservoir may be driven repeatedly, which reveals the same result in theory but poses different experimental challenges~\cite{Oliver:15}.
For a scalar reservoir, or for a scalar observation of the reservoir, consistency is then measured by the \emph{consistency correlation} $\gamma^2$, which is the Pearson-correlation coefficient between the two responses.
For the multivariate response of a large reservoir one may define consistency correlation for each node
\begin{equation}
\gamma_i^2=\langle \bar{x}_i(t)\bar{x}_i'(t)\rangle_t
\label{eq:gammai}
\end{equation}
where $\langle\cdot\rangle_t$ is the average with respect to time $t$, and $\bar{x}$ indicates normalization of $x$ to zero mean and unit variance.
An average over all nodes
\begin{equation}
\hat{\gamma}^2 = \langle \gamma_i^2\rangle_i
\label{eq:gammahat}
\end{equation}
then accounts for the `global' consistency.
It is worth noting that a different replica test could also be designed for each node in the network separately, or for an arbitrary group of nodes as a subset of the whole reservoir.
This would define a different consistency measure which depends on the selection of the subset.
Such a measure would for instance allow to locate the source of inconsistency.
Nevertheless, we will focus in this work only on the replica test for the whole reservoir as outlined before and illustrated in Fig.~\ref{fig:replica}.

The supervised learning procedure in RC creates a set of nodes in a separate readout layer.
The readout signal is typically a linear superposition of the reservoir nodes
\begin{equation}
y(t) = \sum_{i=1}^N R_i x_i(t)+R_0\;.
\label{eq:yt}
\end{equation}
The $R_i$ are components of the readout vector $\mathbf{R}\in\mathbb{R}^N$ which is typically obtained by ridge regression with respect to a target signal $z(t)$.
For the sake of simplicity, we will omit the bias term $R_0$ in the following discussions.
A separate consistency measure can be determined for these readout nodes by
\begin{equation}
\Gamma_R^2=\langle\bar{y}(t)\bar{y}'(t)\rangle_t
\label{eq:GammaR}
\end{equation}
where $y'(t)$ is the readout with the same vector $\mathbf{R}$ applied to the replica reservoir, and normalization is applied in both $\bar{y}(t)$ and $\bar{y}'(t)$.
The readout can be considered an emergent quantity which due to the training adjustments is a special projection of the reservoir dynamics.
Its consistency level thus plays a distinguished role as compared to the individual $\gamma_i^2$.

A broader notion of consistency can be found throughout the community in which a system is said to be consistent if similar inputs lead to similar outputs.
However, this idea overlaps with the approximation property~\cite{Maass:02}, and measures of similarity are not yet specified.
We restrict our investigation to only the output similarity given exact repetitions of the input, as described above, in order to probe for the degree of functional dependency.
This way, the consistency property is distinguished from the approximation property.
Generalizing from the consistency correlation $\gamma^2$, a reasonable measure of similarity among the responses to repeated inputs is given by the consistency correlation $\Gamma^2$ of any \emph{observable} $y(t)=h(\mathbf{x}(t))$ of the system.
The typical reservoir readouts $\mathbf{R}$ are the special case in which the projection $h:\mathbb{R}^N\mapsto\mathbb{R}$ is linear.
This notion is distinct from consistency in readouts which are a filtered function of the reservoir state, for instance $\dot{y}(t)=g(y(t),\mathbf{x}(t))$.
Future work may be oriented towards a general consistency concept including filtered signals, to account for phenomena like rate coding in neuronal circuits where the timing of individual spikes is inconsistent with respect to certain reference signals.
In neuroscience, consistency on the level of spike timing is known as \emph{reliability}~\cite{Mainen:95,Goldobin:06}.
Despite the similarity between consistency and reliability, however, the two concepts are not identical due to the different context and measure of functional dependency.
In nonlinear science, consistency can also be compared to synchronization due to common drive, where noise is often chosen as a driving signal~\cite{Teramae:04,Goldobin:05,Pimenova:16}.
What distinguishes consistency from synchronization phenomena is that it is a property of a single system subject to a driving signal.
Moreover, the consistency property is inherent to the system even if the signal is not repeatedly presented.

We illustrate our consistency concept at the example of an ESN driven by noise.
This reservoir updates in discrete time $t\in\mathbb{Z}$ and reads
\begin{equation}
\label{equation:update}
\mathbf{x}(t+1) = \tanh(W\cdot\mathbf{x}(t) + V\cdot\mathbf{u}(t+1) + \boldsymbol{\beta})
\end{equation}
where $\tanh(\cdot)$ is applied element-wise.
The internal connectivity is summarized in the matrix $W\in\mathbb{R}^{N\times N}$, and the input injection is contained in $V\in\mathbb{R}^{N\times L}$.
We create a network with $N=200$ nodes and connect nodes randomly with a probability of $p=2.5\%$.
The weight of each connection is then chosen from a normal distribution $\mathcal{N}(0,1)$, and the weight is zero if there is no connection.
The resulting matrix $W$ is then scaled by global factor to achieve a desired spectral radius $\rho$.
The spectral radius is a key parameter in the design of ESN which can be thought of as the internal gain of the dynamical system.
The input connections in $V$ are created in a way that the input is injected into each node, with the weight for each connection taken from a uniform distribution between $-1$ and $1$.
$\boldsymbol{\beta}$ is a vector of biases which shifts the operating window of each node to different regions of the $\tanh(\cdot)$-nonlinearity. 
The bias for each node is set to one here, meaning $\boldsymbol{\beta}=(1,1,...,1)^\top$.
The input is chosen to be a scalar ($L=1$) IID random variable taken from a normal distribution with zero mean and unit variance, \(u(t) \sim\mathcal{N}(0,1)\).

\begin{figure}
\begin{center}
\includegraphics[width=\columnwidth]{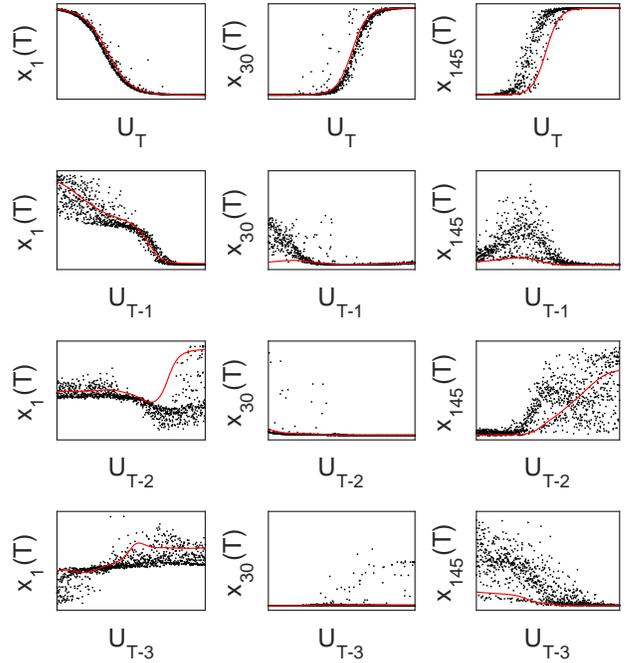}
\end{center}
\caption{Selected sections through the input-output dependency of an ESN generated via repeated drives with a perturbation placed at time $T-\tau$, see text.
From the $N=200$ nodes, activations of three selected nodes are plotted for lags $\tau \in [0,3]$. Black points: Inconsistent regime at $\rho=3$. Red curves: Consistent regime at $\rho=2.2$ in the same network.}
\label{fig:functs}
\end{figure}

In order to visualize sections through the input-output relationship, i.e. between the drive $u(t)$ and the response $\mathbf{x}(t)$, we create first a sufficiently long reference sequence $U(t)=(U_1,...,U_T)$ as a single realization of the noise process.
The ESN is then repeatedly driven with a variation of this sequence, in which only a single element $U_{T-\tau}$ takes a different value in each run.
We select a few nodes to plot the dependency of $x_i(T)$ on this variable for different lags $\tau$.
Figure~\ref{fig:functs} shows these sections for selected nodes in a regime of low consistency at $\rho=3$ together with the same sections in the same network and drive, but at $\rho=2.2$ where the response is completely consistent.
The consistent response reveals slices through a functional $\mathbf{x}(T)=\mathcal{F}[\mathbf{U}](T)$.
In the low-consistency regime, however, the dependency is blurred to an extent which depends on the individual node as well as on the state of the drive.
In both cases, we recognize the $\tanh(\cdot)$-nonlinearity in the instantaneous response, and higher nonlinear features for increasing lag.
Inferring to an exact repetition of the drive, the variability in the response can be seen by taking a vertical slice of the response portraits shown in Fig.~\ref{fig:functs}.
In the consistent case the slice reduces to a point, whereas in the inconsistent case this is a distribution, describing the values the response may take for a certain input sequence.
The consistency correlation thus measures the degree to which the total variability through time arises from a functional dependence as compared to the chaotic variability.


\section{Memory and Consistency}
\label{sec:memory}

Recurrent neural networks allow the input signal to propagate through the network for multiple timesteps, meaning that the current state of the network contains information about the history of the input.
This ability to store past information is a key component of RNN which enables them to be powerful tools for computation on sequential data.
For ESN, understanding the relationship between the memory profile and the hyper-parameters of the network is important for optimal, task specific reservoir design.
In this section we investigate the memory of ESN in the context of consistency.

\begin{figure}
\centering
\includegraphics[width=0.48\textwidth]{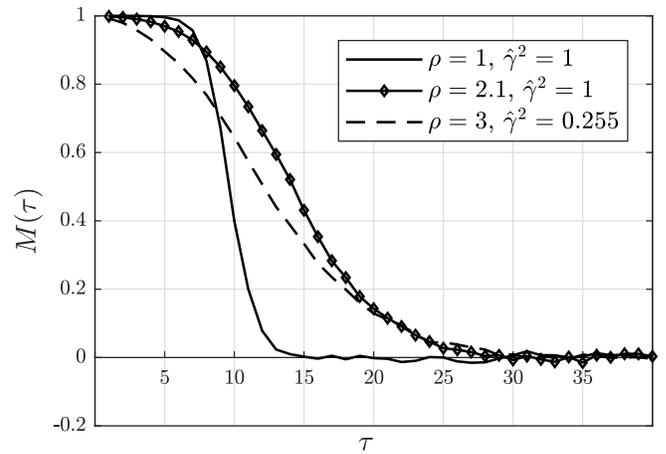}
\caption{\label{fig:1}The memory profile of an ESN in the consistent regime, the transition to inconsistency and in the highly inconsistent regime.
The correlation between the true and predicted lagged input is plotted against the lag.}
\end{figure}

The linear reconstruction task and the associated memory capacity (MC) measure as introduced by Jaeger~\cite{Jaeger:01} are commonly used to quantify the fading memory.
The ESN is trained to reconstruct the input $\tau$ timesteps ago, meaning that the training target is $z_\tau(t)=u(t-\tau)$, $\tau\in\mathbb{N}$.
The reconstruction accuracy at lag $\tau$ is an indicator of the amount of information held in the network about the input at that lag.
With the reconstruction from the ESN reading $y_\tau(t)$, the accuracy is measured by the correlation coefficient
\begin{equation}
M(\tau)=\langle\bar{z}_{\tau}(t)\bar{y}_{\tau}(t)\rangle_t
\label{eq:mtau}
\end{equation}
where the usual normalizations apply.
This performance measure as a function of $\tau$ is the memory profile of the ESN, as shown in Fig.~\ref{fig:1}.
The memory capacity $I_{MC}$ is an integral over all lags which measures the total linear memory of the reservoir~\cite{Jaeger:01}
\begin{equation}
I_{MC} = \sum_{\tau=1}^{\infty}M(\tau)^2\text{ .}
\end{equation}
A key result with regard to ESN memory is that the MC is maximized at the \textit{edge of chaos}, where the maximal Lyapunov exponent becomes positive and the reservoir dynamics change from a stable to an unstable regime.
This has informed part of the design strategy for ESN, which is to scale the spectral radius of the internal weight matrix to just before this point in order to maximize memory.

\begin{figure}
\centering
\includegraphics[width=0.48\textwidth]{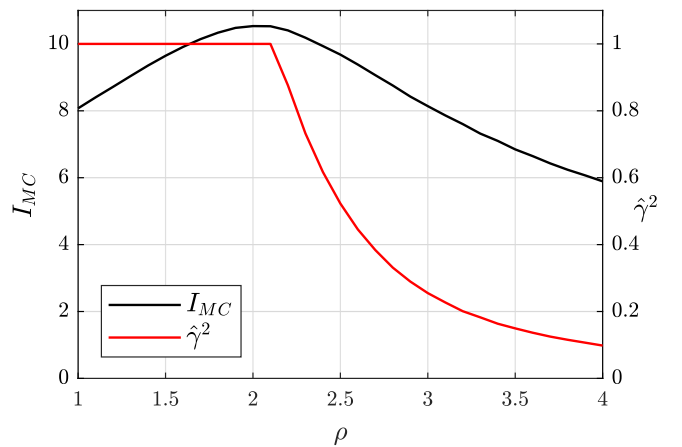}
\caption{\label{fig:2}The memory capacity and consistency as a function of the spectral radius of the internal weight matrix $W$.}
\end{figure}

We investigate here the memory of an ESN in the stable regime, at the transition to instability, and in the unstable regime.
The ESN is set up as in Sec.~\ref{sec:consistency}, but with a size of $N=500$ nodes and wiring probability $p=10\%$.
We perform the reconstruction task for a spectral radius of \(\rho \in [1, 4]\), which encompasses the aforementioned regimes with different levels of consistency.
The uncorrelated input ensures that any memory observed is purely from transient activation in the reservoir, rather than due to any autocorrelation already present in the input.
We measure the global consistency $\hat{\gamma}^2$ according to Eq.~\eqref{eq:gammahat}.
The values are averaged over 10 different realizations after decay of transients.
The results of this experiment are shown in Fig.~\ref{fig:1} and Fig.~\ref{fig:2}. 

Considering consistency in the terminology of ESN, complete consistency (\(\hat{\gamma}^2=1\)) is equivalent to the echo-state property, and the transition to inconsistency at \(\rho \approx 2\) is equivalent to the edge of chaos.
The results in Fig.~\ref{fig:2} support the notion that the memory capacity of an ESN is maximized approximately at the onset of inconsistency.
Starting from no connectivity ($\rho=0$), memory increases with increasing $\rho$ until the network transitions into an unstable regime, leading to a decrease in consistency which leads to a decrease in memory.
However, even deep in the inconsistent regime the network still performs relatively well, with the accuracy of the reproduction for large lags (\(\tau>10\)) being greater than the network in the fully consistent regime (Fig.~\ref{fig:1}).
For the memory capacity measure, the highly inconsistent network performs comparably to the reservoir scaled to \(\rho = 1\) (Fig.~\ref{fig:2}), which was the previous best practice in reservoir design.
This is a surprising result, as one may expect that a significant loss of consistency should be accompanied by a comparable overall loss of reconstruction accuracy.
Particularly for large lags, the effect of inconsistency is expected to accumulate, because the signal-to-noise ratio is effectively reduced with every propagation step.
However, it is at large lags that the performance of the inconsistent reservoir is the strongest relative to the consistent and transition cases.
Thus it seems wrong to assume that inconsistency simply acts like additive noise at each node.
There must instead be a mechanism which enables the input to propagate through an inconsistent reservoir in a way that is recoverable by a linear readout.

\begin{figure}
\centering
\includegraphics[width=0.48\textwidth]{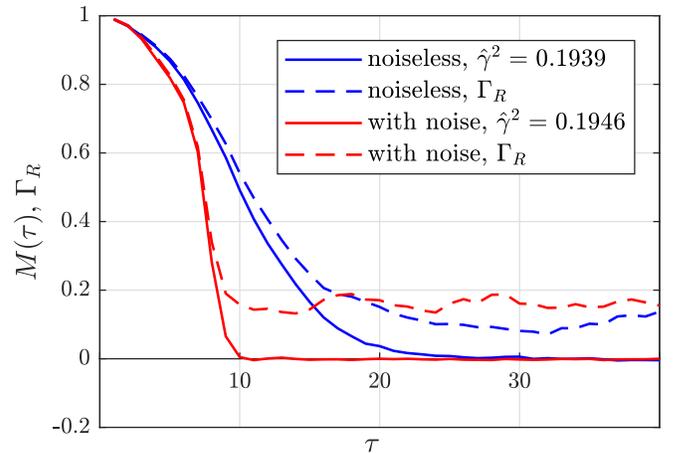}
\caption{\label{fig:3}The memory profile (solid lines) of an ESN (Eq.~\eqref{equation:update}) which is in the inconsistent regime due to large spectral radius (\(\rho = 3\), blue) and an ESN which has been brought into the inconsistent regime by the introduction of noise as per Eq.~\eqref{equation:noisey_ESN} (\(\rho = 1\), red).
The network parameters were chosen so that the two systems have comparable global consistencies (0.1939 and 0.1946).
The dashed lines show the square root of the output-consistency correlation, $\Gamma_R$.}
\end{figure}

We obtain further insight into the modes of signal propagation by a different numerical experiment.
Starting from an ESN which is initially in the fully consistent regime, say \(\rho = 1\), we induce inconsistency by introducing a source of noise at each node.
The update equation for the reservoir becomes
\begin{equation}
\label{equation:noisey_ESN}
\begin{split}
\mathbf{x}(t+1) = \tanh((1-r)(W\cdot \mathbf{x}(t) + V\cdot u(t) + \boldsymbol{\beta}) \\ +r \boldsymbol{\xi}(t))
\end{split}
\end{equation}
where \(r\) is a parameter that determines the amount of noise, and $\xi_i(t)\sim\mathcal{N}(0,1)$.
Figure~\ref{fig:3} shows the memory profile of two ESN, each in an inconsistent regime with similar levels of global consistency (\(\hat{\gamma}^2\approx 0.2\)).
One ESN has inconsistency naturally due to instability in the autonomous reservoir dynamics caused by a large spectral radius.
The other has inconsistency introduced via noise.
The results show that on the memory test the standard reservoir in an inconsistent regime performs better than a reservoir spiked with noise.
Figure~\ref{fig:3} also shows the square root of the output consistency, $\Gamma_R$, according to Eq.~\eqref{eq:GammaR} with $\mathbf{R}$ depending on $\tau$.
For the outputs associated with reproducing recent inputs, the output consistency is much larger than the consistency of the network as a whole.
This means that there are readout projections $\mathbf{R}$ for which the response of the ESN is highly consistent, even if the network as a whole is inconsistent.
We also see that the accuracy of the memory reconstruction is closely bounded by the $\Gamma_R$ values.
This is in agreement with consistency theory~\cite{Jungling:18}, meaning that the trained readout exploits consistency optimally.

\section{Consistency Profile}
\label{sec:profile}

\begin{figure}
\centering
\includegraphics[width=0.48\textwidth]{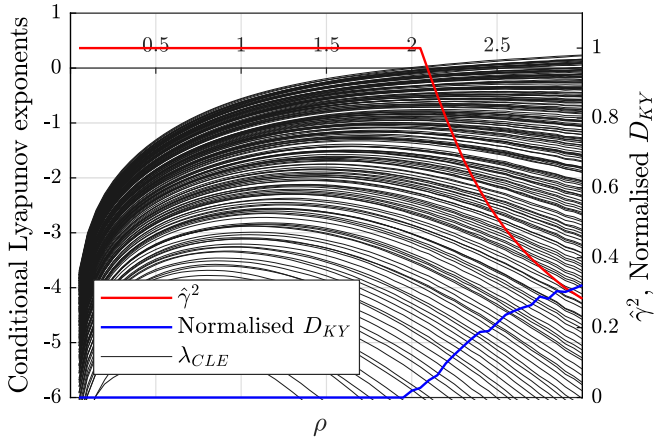}
\caption{\label{fig:4}The Lyapunov spectrum (left axis, black) of an ESN with $N=200$ nodes plotted against the spectral radius of the network. The consistency (right axis, red) and normalised Kaplan-Yorke dimension (right axis, blue) are also plotted. We see that the transition into inconsistency occurs as the maximal Lyapunov exponent becomes positive. However, even for high inconsistency (\(\rho=3,~\hat{\gamma}^2 \approx 0.3\)) approximately $90\%$ of Lyapunov exponents remain negative. 
}
\end{figure}

The disturbance of the propagating signal due to chaos emerges to be less dramatic than by a source of noise at every node.
This is little surprising when we take into account that chaos effectively populates only a few degrees of freedom according to the attractor dimension.
We follow this idea by first calculating the conditional Lyapunov spectrum for an ESN~\cite{Verstraeten:09}.
The transfer of the dimension to a driven system is possible if we interpret chaotic dimensionality as additional degrees of freedom superposed to the signal response.
Figure~\ref{fig:4} shows the CLE calculated via the Gram-Schmidt procedure~\cite{eckmann1985ergodic} together with the global consistency $\hat{\gamma}^2$ and the conditional attractor dimension $D_{KY}$, which is the Kaplan-Yorke dimension from the CLE~\cite{frederickson1983liapunov}.
As expected, the transition to inconsistency corresponds to the crossing of the maximal Lyapunov exponent from negative to positive.
However, even for large inconsistency approximately $90\%$ of the Lyapunov exponents remain negative, and the chaotic dimension is still small compared to the state-space dimension $N$.
Thus the ESN is still effectively stable in a large portion of the available directions.
These more stable directions may have a higher level of consistency than others, which gives a first idea of why a signal can still propagate through a globally inconsistent network for many timesteps.

The Lyapunov dimension provides only a very basic argument on the distribution of signal and noise in the response of the reservoir.
In general, Lyapunov exponents are little related to correlations in dynamical systems.
This is because the attractor dimension is a topological property, whereas correlations are \emph{geometrical} properties of the dynamics.
The chaotic degrees of freedom of the reservoir are not confined to a trivial subspace, as for instance in the case that the reservoir was a linear system, e.g. by omitting the $\tanh(\cdot)$-nonlinearity.
Lyapunov vectors from the nonlinear Eq.~\eqref{equation:update} are time dependent and effectively distribute the chaotic instabilities over all degrees of freedom.
In the following, we introduce a comprehensive characterization of the distribution of signal and chaos (noise), including the effect of regularization, based on principal-component analysis (PCA).

\begin{figure}
	\centering
	\includegraphics[width=0.48\textwidth]{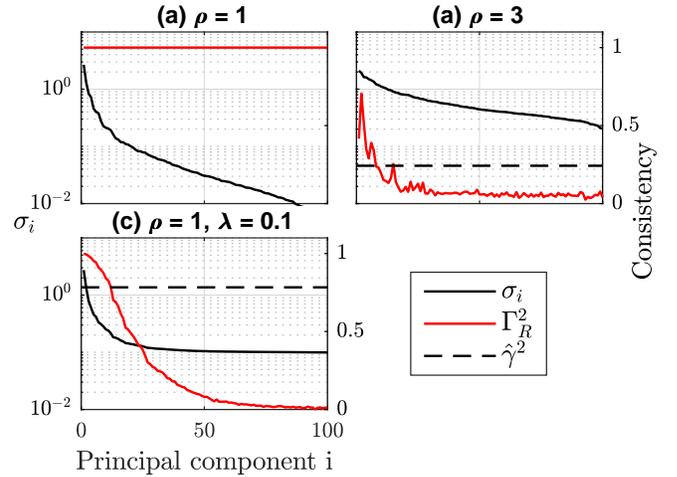}
	\caption{\label{fig:6}The response profile and principal component readout consistencies for an ESN in the consistent regime (top left), in the inconsistent regime (top right), and in the consistent regime with regularization (bottom left).}
\end{figure}

We first apply the PCA directly to the full response of the ESN under various conditions.
This is done by performing singular-value decomposition on the covariance matrix for $\mathbf{x}(t)$,
\begin{equation}
\begin{split}
C_{xx} &= \begin{bmatrix}
    \langle x_1 x_1\rangle & \langle x_1 x_2\rangle & \dots  & \langle x_1 x_N\rangle \\
    \langle x_2 x_1\rangle & \langle x_2 x_2\rangle & \dots  & \langle x_2 x_N\rangle \\
    \vdots & \vdots & \ddots & \vdots \\
    \langle x_N x_1\rangle & \langle x_N x_2\rangle & \dots  & \langle x_N x_N\rangle
    \end{bmatrix}\\
&= Q\Sigma^2 Q^\top\text{ .}
\end{split}
\label{eq:corr}
\end{equation}
Note that in contrast to the correlation functions before, here we do not apply normalization.
The columns $\mathbf{Q}_i$ of $Q$ form an orthonormal set of vectors in the direction of the principal components (PCs) of the response.
The diagonal matrix $\Sigma$ contains the sizes $\sigma_i$ of the principal components, which measure the extent of the response of $\mathbf{x}(t)$ in the corresponding PC-direction.
Figure~\ref{fig:6} shows the PC profiles for different ESN responses.
Moreover, by taking the PC directions as readouts, $\mathbf{R}_i=\mathbf{Q}_i^\top$, we calculate the corresponding readout consistency correlations $\Gamma_i^2$ for these directions and thus obtain a consistency distribution.
We apply this procedure to a network in the completely consistent ($\rho =1 $) and in the inconsistent ($\rho=3$) regime (Fig.~\ref{fig:6}a-b).
We further consider the effect of regularization, by adding measurement noise to the reservoir state, $\mathbf{x}(t)\rightarrow\mathbf{x}(t) + \lambda\boldsymbol{\xi}(t)$, $\xi_i(t)\sim\mathcal{N}(0,1)$, equivalent to performing ridge regression (Fig.~\ref{fig:6}c).

For the completely consistent case we find $\Gamma_i^2$ = 1 for each readout, as expected.
The effect of regularization is to remove the smallest $\sigma_i$ of the response, leaving only those directions active in which the added noise is small relative to the propagated input signal.
The resulting consistency distribution reveals a new way to describe the capacity of the ESN, which turns out to be significantly smaller than the number of degrees of freedom of the reservoir.
This is important as ridge regression is commonly used when training ESN, and in physical cases of reservoir computing there is often measurement noise which acts in a similar way to regularization.
For the inconsistent case, the directions of larger response tend to have consistency above the global consistency $\hat{\gamma}^2$, and vice versa. 
However, even in the direction with greatest consistency, the response does not reach the consistency level observed for readouts in the memory profile (Fig.~\ref{fig:4}).

We describe in the following how to trace directions of particular consistency levels.
The notion of consistency directly leads to a new characteristic set of readouts, different from the orthogonal set obtained from the PCA.
We perform an alternative PCA on the consistent component of the reservoir response~\cite{Jungling:18}.
The consistent component can be found via an average over an ensemble of replica states $\mathbf{x}^c(t) = \langle \mathbf{x}^i(t) \rangle_i$.
Following naturally from this, an inconsistent component $\mathbf{n}^i(t)$ can be defined for each replica such that $\mathbf{x}^i(t) = \mathbf{x}^c(t) + \mathbf{n}^i(t)$.
The only underlying assumption in this decomposition is ergodicity, which we found is to a high degree satisfied in large reservoirs.
Thus $\mathbf{n}^i(t)$ becomes equivalent to a realization of a noise-like process which, for long enough time, does not correlate with either the consistent component or other realizations of this process~\cite{Jungling:18}. 
This allows us to relate the variance of the consistent component to the covariance of two replica states (Eq.~\eqref{equation:consistent_component_correlation}).
\begin{align}
\langle \mathbf{x}^i\mathbf{x}^j\rangle &= \langle (\mathbf{x}^c+\mathbf{n}^i)(\mathbf{x}^c+\mathbf{n}^j)\rangle \nonumber \\
\langle \mathbf{x}^i\mathbf{x}^j\rangle &= \langle \mathbf{x}^c\mathbf{x}^c\rangle + \langle \mathbf{x}^c\mathbf{n}^j\rangle +\langle \mathbf{n}^i\mathbf{x}^c\rangle + \langle \mathbf{n}^i\mathbf{n}^j\rangle \label{equation:consistent_component_correlation} \\
\langle \mathbf{x}^i\mathbf{x}^j\rangle &= \langle \mathbf{x}^c\mathbf{x}^c\rangle \nonumber
\end{align}

This result extends to the covariance matrix of the consistent component, which is equal to the cross-covariance matrix between two replicas in the long time limit.
Thus, through PCA of the cross-covariance matrix, we can access the principal components of the consistent part of the response, which then can be compared with the full response of the ESN.
The cross-covariance matrix of two replica responses $\mathbf{x}(t)$ and $\mathbf{x}'(t)$ reads
\begin{equation}
\begin{split}
C_{xx'} &= \begin{bmatrix}
    \langle x_1 x'_1\rangle & \langle x_1 x'_2\rangle & \dots  & \langle x_1 x'_N\rangle \\
    \langle x_2 x'_1\rangle & \langle x_2 x'_2\rangle & \dots  & \langle x_2 x'_N\rangle \\
    \vdots & \vdots & \ddots & \vdots \\
    \langle x_N x'_1\rangle & \langle x_N x'_2\rangle & \dots  & \langle x_N x'_N\rangle
    \end{bmatrix}\\
    &=C_{c}
\end{split}
\label{eq:corr2}
\end{equation}
where $C_c$ is the covariance matrix of $\mathbf{x}^c(t)$, and no normalizations are applied.
Our numerical experiments confirm that the ergodicity assumption behind this equality holds well, and finite-size effects are negligible with reasonable time-series lengths.
We perform SVD on both $C_{xx}$ (Eq.~\eqref{eq:corr}) and $C_c$ (Eq.~\eqref{eq:corr2}) to find the principal components, using the symmetry and positive definiteness of the covariance matrices
\begin{align}
    C_{xx} &= Q_{xx} \Sigma_{xx}^2 Q_{xx}^\top \\
    C_c &= Q_c \Sigma_c^2 Q_c^\top
\end{align}
From here, we aim to measure the relative orientation of the signal response within the full response, which can be geometrically illustrated as two nested ellipsoids.

\begin{figure}
\centering
\includegraphics[width=0.48\textwidth]{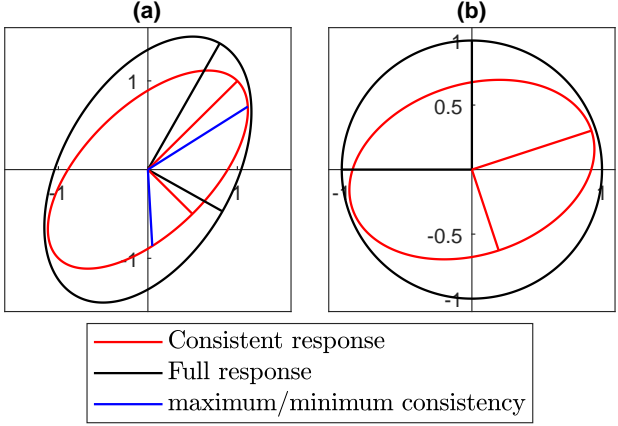}
\caption{\label{fig:5}Component ellipses of the full response (black) and the consistent response (red) for the test system Eq.~\eqref{test_system} with directions of maximum and minimum consistency (blue).
(a) Original coordinates $\mathbf{x}(t)$.
(b) Coordinates $\bar{\mathbf{x}}(t)$ normalized as per Eq.~\eqref{eq:To} to an isotropic full response.}
\end{figure}

To demonstrate this method we will consider a simple two dimensional test system comprised of a consistent and inconsistent component.
\begin{align}
\mathbf{x}^c(t) &= \xi_1(t)\begin{bmatrix} 1 \\
           1 
        \end{bmatrix} + \xi_2(t)\begin{bmatrix} 0.5 \\
           -0.5
        \end{bmatrix} \nonumber \\
        \mathbf{n}^i(t) &= \nu_{1,i}(t)\begin{bmatrix} 1 \\
           0 
        \end{bmatrix} + \nu_{2,i}(t)\begin{bmatrix} 0 \\
           0.3 
        \end{bmatrix} \label{test_system}\\
        \mathbf{x}^1(t) &= \mathbf{x}^c(t) + \mathbf{n}^1(t) \nonumber \\
        \mathbf{x}^2(t) &= \mathbf{x}^c(t) + \mathbf{n}^2(t) \nonumber
\end{align}
All $\xi_j(t),\nu_{j,i}(t)\sim\mathcal{N}(0,1)$.
The results of PCA on this test system are shown in Fig.~\ref{fig:5}a, where the components span ellipses.
Besides the main axes of the full and consistent responses, we also show the directions of maximum and minimum consistency.
These are the directions in which the ratio between inner and outer ellipse are maximal or minimal, respectively.
The consistency directions do not align with the orthogonal PC axes of either ellipse.
In order to define the consistency directions, we introduce a coordinate transformation which normalizes the full response components
\begin{equation}
T_{\circ} = Q_{xx} \Sigma_{xx}^{-1} Q_{xx}^\top\;.
\label{eq:To}
\end{equation}
We apply $T_{\circ}$ to the two replica states $\mathbf{x}^1(t)$ and $\mathbf{x}^2(t)$ to get $\mathbf{\bar{x}}^1(t)=T_{\circ} \mathbf{x}^1(t)$ and $\mathbf{\bar{x}}^2(t)=T_{\circ} \mathbf{x}^2(t)$.
The bar notation indicates the transformed states.
This transformation preserves the relative proportions relevant for consistency.
The normalized geometry is shown in Fig.~\ref{fig:5}b.
The full response component ellipse has been transformed into a unit circle.
The directions of maximum and minimum consistency here align with the principal components of the consistent response, and moreover, also with the components of the inconsistent response (not shown).
In other words, the result of the isotropic reservoir response in the new coordinates is the consistent and the inconsistent part being complementary with a shared set of principal axes. 
Moreover, the diagonal elements of $\overline{\Sigma}_c^2$ in these coordinates are the consistency levels in the corresponding directions.

\begin{figure}
\centering
\includegraphics[width=0.48\textwidth]{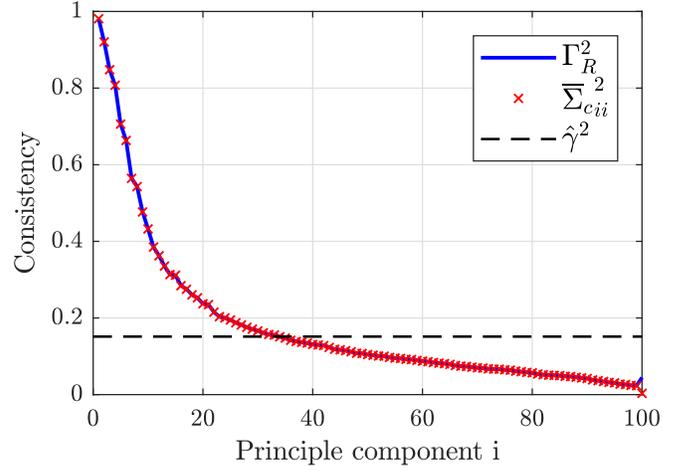}
\caption{\label{fig:7}Consistency profile of ESN with $N=100$ and $\hat{\gamma}^2=0.15$. The readout consistency along the principal components (blue) and the consistency taken from $\overline{\Sigma}_c$ (red) agree strongly. }
\end{figure}

We apply this method to an ESN with 100 nodes and a global consistency of $\hat{\gamma}^2 = 0.15$.
By plotting the consistency along each principal component, we get the consistency profile (Fig.~\ref{fig:7}).
In approximately $30\%$ of directions the consistency is larger than the global consistency of the ESN, with some directions maintaining a very high level of consistency.
This suggests that it is through these directions of higher consistency that the input signal is able to propagate in a manner that is recoverable by the linear readout. 


\section{Conclusion}

We have applied the concept of consistency to echo-state networks as a new way of characterizing the echo-state property. 
In particular, we have assessed the performance of ESN in a range of consistency regimes on a memory task.
We found that inconsistency is not as destructive to fading memory as expected, with the inconsistent reservoir performing comparably, and even outperforming consistent reservoirs.
The reason for signals surviving in the inconsistent regime is found in the distribution of signal and noise. 
We introduced the consistency profile based on principal-component analysis as a portrait of the high-dimensional response.
We found that a few directions of high consistency are always present, bypassing the chaotic or noise-induced fluctuations.
Our method is applicable to an arbitrary reservoir computer, including physical media in which noise is inherently present.

While complete consistency is always desirable for reservoir design, we found that in typical ESN with regularization the response is poorly exploited.
Our findings thus may give rise to unsupervised pre-training methods which aim to maximize the consistent dimension in a response subject to chaos, noise, and regularization.
Furthermore, one may think of shared processing of multiple input channels, in which the distinct inputs act as sources of inconsistency for each other rather than complementary channels to synthesize the desired output.
This is likely to be the case in many biological examples, where a single computational unit receives multiple stimuli.
Our method may help to understand how computational capacity is distributed and routed in such systems.
In summary, the consistency profile including effects of regularization may lead to an enhanced understanding of computational capacity in noisy neuronal microcircuits, and also prove useful in unsupervised optimization procedures for reservoir design.


\begin{acknowledgments}
TL is supported by the Australian Government Research Training Program at The University of Western Australia.
AK was supported by the Hackett Postgraduate Scholarship.
MS is supported by Australian Research Council Discovery Project DP180100718.
\end{acknowledgments}

\section*{References}
%


\end{document}